\documentclass[12pt]{article}

\usepackage{amsthm,amsmath,amssymb,url,cite}
\newtheorem{thm}{Theorem}[section]
\newtheorem{lem}[thm]{Lemma}
\newtheorem{cor}[thm]{Corollary}
\newtheorem{prop}[thm]{Proposition}
\theoremstyle{remark}
\newtheorem*{rmk}{Remark}

\newtheorem*{eg*}{Example}

\numberwithin{equation}{section}

\renewcommand{\qed}{{\hfill\rule{4pt}{7pt}}\medskip}

\def\pf{\noindent {\it Proof.} }

\begin{document}

\begin{center}
{\Large\bf Short Proofs of Summation and Transformation\\[5pt]
 Formulas for Basic Hypergeometric Series}\footnote{to appear in 
 {\it J. Math. Anal. Appl.}}
\end{center}
\vskip 2mm
\centerline{Victor J. W. Guo$^1$  and Jiang Zeng$^2$}

\begin{center}
Institut Camille Jordan,
Universit\'e Claude Bernard (Lyon I)\\
F-69622, Villeurbanne Cedex, France \\
{\tt $^1$guo@math.univ-lyon1.fr, http://math.univ-lyon1.fr/\textasciitilde{guo}}\\
{\tt $^2$zeng@math.univ-lyon1.fr, http://math.univ-lyon1.fr/\textasciitilde{zeng}}
\end{center}

\vskip 0.7cm {\small \noindent{\bf Abstract.}
We show that several terminating summation and transformation formulas for
basic hypergeometric series can be proved in a straightforward way.
Along the same line, new finite forms of Jacobi's triple product
identity and Watson's quintuple product identity are also proved.
}


%\vskip 0.2mm
%\noindent{\bf MR Subject Classifications}: Primary 33D15; %%Secondary 05A30

\section{Introduction}
We follow the standard notation for $q$-series and basic
hypergeometric series in
\cite{GR}. The {\it $q$-shifted factorial} is defined by
\[
(a;q)_\infty=\prod_{k=0}^{\infty}(1-aq^{k}),\quad\text{and}
\quad (a;q)_n=\frac{(a;q)_\infty}{(aq^n;q)_\infty}, \quad\text{for}
\quad n\in\mathbb{Z}.
\]
As usual, we employ the abbreviated notation
$$
(a_1,a_2,\ldots,a_m;q)_n
=(a_1;q)_n (a_2;q)_n \cdots(a_m;q)_n,
\quad\text{for}\quad n=\infty\quad\text{or}\quad n\in\mathbb{Z}.
$$
The {\it basic hypergeometric series $_{r+1}\phi_r$} is defined as
$$
_{r+1}\phi_{r}\left[\begin{array}{c}
a_1,a_2,\ldots,a_{r+1}\\
b_1,b_2,\ldots,b_{r}
\end{array};q,\, z
\right]
=\sum_{k=0}^{\infty}\frac{(a_1,a_2,\ldots,a_{r+1};q)_k z^k}
{(q,b_1,b_2,\ldots,b_{r};q)_k}.
$$
An $_{r+1}\phi_{r}$ series is called {\it well-poised} if
$a_1q=a_2b_1=\cdots=a_{r+1}b_r$ and {\it very-well-poised}
if it is well-poised and $a_2=-a_3=q\sqrt{a_1}$.

The starting point of this paper is the observation that
the  $k$th term of a well-poised hypergeometric series
\begin{align*}
F_k(a_1,a_2,\ldots,a_{r+1};q,z):=\frac{(a_1,\ldots,a_{r+1};q)_k}
{(q,a_1q/a_2,\ldots,a_1q/a_{r+1};q)_k} z^k
\end{align*}
satisfies the following contiguous relations:
\begin{align}
F_k(a_1,a_2,\ldots,a_rq, a_{r+1};q,z)&-F_k(a_1,a_2,\ldots,a_r, a_{r+1}q;q,z)\nonumber \\
&=\alpha F_{k-1}(a_1q^2,a_2q,\ldots,a_{r+1}q;q,z),
\label{eq:contiguous1} \\[5pt]
F_k(a_1,a_2,\ldots,a_r, a_{r+1};q,qz)&-F_k(a_1,a_2,\ldots,a_r, a_{r+1}q;q,z)\nonumber\\
&=\beta F_{k-1}(a_1q^2,a_2q,\ldots,a_{r+1}q;q,z), \label{eq:contiguous2}
\end{align}
where
\begin{align*}
\alpha &=\frac{(a_r-a_{r+1})(1-a_1/a_ra_{r+1})(1-a_1)(1-a_1q)(1-a_2)
\cdots (1-a_{r-1})z}{(1-a_1/a_r)(1-a_1/a_{r+1})(1-a_1q/a_2)
\cdots (1-a_1q/a_{r+1})}, \\[5pt]
\beta &=-\frac{(1-a_1)(1-a_1q)(1-a_2)\cdots (1-a_r)z}
{(1-a_1/a_{r+1})(1-a_1q/a_2)\cdots (1-a_1q/a_{r+1})}.
\end{align*}
Indeed, \eqref{eq:contiguous1} and \eqref{eq:contiguous2} correspond
 respectively to the following trivial identities:
\begin{align*}
%\frac{(1-bx)(1-c)}{(1-a/b)(1-ax/c)}-\frac{(1-b)(1-cx)}{(1-ax/b)(1-a/c)}
%&=\frac{(b-c)(1-a/bc)(1-ax)(1-x)}{(1-a/b)(1-a/c)(1-ax/b)(1-ax/c)}\\[5pt]
\frac{(1-bx)(1-ax/b)}{(1-b)(1-a/b)}-\frac{(1-cx)(1-ax/c)}{(1-c)(1-a/c)}
&=\frac{(b-c)(1-a/bc)(1-x)(1-ax)}{(1-b)(1-c)(1-a/b)(1-a/c)},\\[5pt]
x-\frac{(1-cx)(1-ax/c)}{(1-c)(1-a/c)}
&=-\frac{(1-x)(1-ax)}{(1-c)(1-a/c)},
\end{align*}
with $a=a_1$, $b=a_r$, $c=a_{r+1}$ and $x=q^k$.

In this paper
we shall give one-line human proofs of several important well-poised
$q$-series identities based on the above contiguous relations.
This is a continuation of our previous works~\cite{Guo,GZ}, where some simpler
$q$-identities were proved by this technique.

Given a summation formula
\begin{align}
\sum_{k=0}^{n}F_{n,k}(a_1,\ldots,a_s)=S_n(a_1,\ldots,a_s),
\label{eq:fnk-sum}
\end{align}
where $F_{n,k}(a_1,\ldots,a_s)=0$ if $k<0$ or $k>n$,
if one can show that the summand $F_{n,k}(a_1,\ldots,a_s)$
satisfies the following recurrence relation:
\begin{align}
F_{n,k}(a_1,\ldots,a_s)-F_{n-1,k}(a_1,\ldots,a_s)
=\gamma_n F_{n-1,k-1}(b_1,\ldots,b_s)
\label{eq:fnk-ind}
\end{align}
for some parameters $b_1,\ldots,b_s$,
where $\gamma_n$ is independent of $k$,
then the proof of the identity \eqref{eq:fnk-sum} is completed by
induction if one can show that $S_n(a_1,\ldots,a_s)$
satisfies the following recurrence relation
\begin{align}
S_{n}(a_1,\ldots,a_s)-S_{n-1}(a_1,\ldots,a_s)=\gamma_n S_{n-1}(b_1,\ldots,b_s). \label{eq:ssn}
\end{align}

If $S_n(a_1,\ldots,a_s)$ appears as a {\it closed form}
as in Jackson's $_{8}\phi_7$ summation formula \eqref{eq:8phi7},
 then the verification of  recursion \eqref{eq:ssn} is  routine.
If $S_n(a_1,\ldots,a_s)$ appears as a sum of $n$ terms, i.e.,
\begin{align}
S_n(a_1,\ldots,a_s)=\sum_{k=0}^{n}G_{n,k}(a_1,\ldots,a_s),
\label{eq:gnk-sum}
\end{align}
where $G_{n,k}(a_1,\ldots,a_s)=0$ if $k<0$ or $k>n$,
we may apply $q$-Gosper's algorithm \cite[p.~75]{Koepf} to
find a sequence $H_{n,k}(a_1,\ldots,a_s)$ of closed forms such that
\begin{align}
%&\hskip -3mm
G_{n,k}(a_1,\ldots,a_s)-G_{n-1,k}(a_1,\ldots,a_s)-\gamma_n G_{n-1,k-1}(b_1,\ldots,b_s)
%\nonumber\\[5pt]
=H_{n,k}-H_{n,k-1}\label{eq:hnk}
\end{align}
and $H_{n,n}=H_{n,-1}=0$, then by telescoping
we get  \eqref{eq:ssn}. Clearly, we have
$$
H_{n,k}=\sum_{j=0}^{k}\left(G_{n,j}(a_1,\ldots,a_s)-G_{n-1,j}(a_1,\ldots,a_s)
-\gamma_n G_{n-1,j-1}(b_1,\ldots,b_s)\right).
$$

In general, we cannot expect \eqref{eq:fnk-ind} or \eqref{eq:hnk} to happen, but as we will show in this paper,
quite a few formulas for basic hypergeometric series can be proved in this way,
such as Jackson's $_{8}\phi_7$
summation, Watson's $_{8}\phi_7$ transformation,
Bailey's $_{10}\phi_9$ transformation, Singh's quadratic transformation,
and a $C_r$ extension of Jackson's $_8\phi_7$ sum due
to Schlosser. The same method can also be used to prove a new finite form of
Jacobi's triple product identity and a finite form of Watson's quintuple product identity.

Since all identities are trivial when $n=0$, we will only indicate the
corresponding recurrence relations \eqref{eq:fnk-ind} (or, in addition, \eqref{eq:hnk})
in our proofs.

It is worth noticing that \eqref{eq:fnk-ind} is not a special case of
Sister Celine's method \cite[p.~58, (4.3.1)]{PWZ}
due to the change of parameters $a_i\to b_i$
 in  the right-hand side.
 For the same reason, Eq.~\eqref{eq:hnk}
is not a hybrid of Zeilberger's method and Sister Celine's method.

%%%%%%%%%%%%%%%%%%%%%%%%%%%%%%%%%%%%%%%%%%%%%%%%%%%%%%%%%%%%%%%%%%%%%%%%%%%%%%%%%%%%%
\section{Jackson's $_{8}\phi_7$ summation formula}\label{sec:Jackson}
Jackson \cite{Jackson21} (see \cite[Appendix (II.22)]{GR}) obtained
%the following $q$-analogue of Dougall's $_7F_6$ sum
a summation formula for a terminating $_{8}\phi_7$ series,
which is one of the most powerful results
in the theory and application of basic hypergeometric series.
\begin{thm}[Jackson's classical $_{8}\phi_7$ summation]For $n\geq 0$, there holds
\begin{align}
 _8\phi_7
\left[\begin{array}{c}
a,\, qa^{\frac12},\, -qa^{\frac12},\, b,\, c,\, d,\, e,\, q^{-n} \\
a^{\frac12},\, -a^{\frac12},\, aq/b,\, aq/c,\, aq/d,\, aq/e,\, aq^{n+1}\end{array};q,\,q
\right]
=\frac{(aq,aq/bc,aq/bd,aq/cd;q)_n}{(aq/b,aq/c,aq/d,aq/bcd;q)_n},
\label{eq:8phi7}
\end{align}
where $a^2q^{n+1}=bcde$.
\end{thm}
\pf Let
$$
F_{n,k}(a,b,c,d,q)=\frac{1-aq^{2k}}{1-a}\frac{(a,b,c,d,a^2q^{n+1}/bcd,q^{-n};q)_k}
{(q,aq/b,aq/c,aq/d,bcdq^{-n}/a,aq^{n+1};q)_k}q^k
$$
be the $k$-th term in Jackson's $_{8}\phi_7$ series. Applying the contiguous relation
\eqref{eq:contiguous1} with $a_r=a^2q^{n}/bcd$ and $a_{r+1}=q^{-n}$, we see that
\begin{align*}
F_{n,k}(a,b,c,d,q)-F_{n-1,k}(a,b,c,d,q)=\alpha_n F_{n-1,k-1}(aq^2,bq,cq,dq,q)
\end{align*}
where
\begin{align*}
\alpha_n=\frac{(a^2q^n/bcd-q^{-n})(1-bcd/a)(1-aq)(1-aq^2)(1-b)(1-c)(1-d)q}
{(1-bcdq^{-n}/a)(1-aq^n)(1-aq/b)(1-aq/c)(1-aq/d)(1-bcdq^{1-n}/a)(1-aq^{n+1})}.
\tag*{\qed}
\end{align*}

\begin{rmk}
A different inductive proof of Jackson's $_{8}\phi_7$ summation formula
can be found in \cite[p.~95]{Slater}.
Substituting $e=a^2q^{n+1}/bcd$ into \eqref{eq:8phi7} and letting $d\to\infty$,
we obtain Jackson's ${}_6\phi_5$ summation (see \cite[Appendix (II.21)]{GR}):
\begin{align}
_6\phi_5
\left[\begin{array}{c}
a,\, qa^{\frac12},\, -qa^{\frac12},\, b,\, c,\, q^{-n} \\
a^{\frac12},\, -a^{\frac12},\, aq/b,\, aq/c,\, aq^{n+1}\end{array};q,\,\frac{aq^{n+1}}{bc}
\right]
=\frac{(aq,aq/bc;q)_n}{(aq/b,aq/c;q)_n}. \label{eq:6phi5}
\end{align}
\end{rmk}
%%%%%%%%%%%%%%%%%%%%%%%%%%%%%%%%%%%%%%%%%%%%%%%%%%%%%%%%%%%%%%%%%%%%%%%%%%%%%%%%%
\section{Watson's $_8\phi_7$ transformation formula}\label{sec:Watson}
Watson's $_8\phi_7$ transformation (see, for example, \cite[Appendix (III.18)]{GR})
formula may be stated as:
\begin{thm}[Watson's classical $q$-Whipple transformation]For $n\geq 0$, there holds
\begin{align}
&\hskip -2mm  _8\phi_7
\left[\begin{array}{c}
a,\, qa^{\frac12},\, -qa^{\frac12},\, b,\, c,\, d,\, e,\, q^{-n} \\
a^{\frac12},\, -a^{\frac12},\, aq/b,\, aq/c,\, aq/d,\, aq/e,\, aq^{n+1}\end{array};
q,\,\frac{a^2q^{n+2}}{bcde}
\right]   \nonumber\\[5pt]
& =\frac{(aq,aq/de;q)_n}{(aq/d,aq/e;q)_n}
 {}_4\phi_3
\left[\begin{array}{c}
aq/bc, \, d,\, e,\, q^{-n} \\
aq/b,\, aq/c,\, deq^{-n}/a,\end{array};q,\,q\right].
\label{eq:watson8phi7}
\end{align}
\end{thm}

\pf Let
\begin{align*}
F_{n,k}(a,b,c,d,e,q)
=F_k(a,qa^{\frac12},-qa^{\frac12},b,c,d,e,q^{-n};q,{a^2q^{n+2}}/{bcde}).
\end{align*}
Applying the contiguous relation \eqref{eq:contiguous2}
with $a_{r+1}=q^{-n}$ and $z=\frac{a^2q^{n+1}}{bcde}$, we have
\begin{align*}
F_{n,k}(a,b,c,d,e,q)-F_{n-1,k}(a,b,c,d,e,q)
=\beta_n F_{n-1,k-1}(aq^2,bq,cq,dq,eq,q),
%\label{eq:rec-ln}
\end{align*}
where
$$
\beta_n=-\frac{(1-aq)(1-aq^2)(1-b)(1-c)(1-d)(1-e)a^2q^{n+1}}
{(1-aq/b)(1-aq/c)(1-aq/d)(1-aq/e)(1-aq^{n})(1-aq^{n+1})bcde}.
$$

Let
$$
G_{n,k}(a,b,c,d,e,q)=\frac{(aq,aq/de;q)_n}{(aq/d,aq/e;q)_n}
\frac{(aq/bc,d,e,q^{-n};q)_k}{(q,aq/b,aq/c,deq^{-n}/a;q)_k}q^k.
$$
Then $G_{n,k}(a,b,c,d,e,q)$ satisfies
\begin{align*}
&\hskip -3mm
G_{n,k}(a,b,c,d,e,q)-G_{n-1,k}(a,b,c,d,e,q)-\beta_nG_{n-1,k-1}(aq^2,bq,cq,dq,eq,q) \\
&=H_{n,k}-H_{n,k-1},
\end{align*}
where
\begin{align*}
H_{n,k}=
\frac{(aq;q)_{n-1}(aq/de;q)_n}{(aq/d,aq/e;q)_n}
\frac{(aq/bc,q^{1-n};q)_k (d,e;q)_{k+1}}
{(q,aq/b,aq/c;q)_k (deq^{-n}/a;q)_{k+1}}. \tag*{\qed}
\end{align*}

\begin{rmk}
The order of our recurrence relation is lower than that generated by $q$-Zeilberger's algorithm or
the $q$-WZ method \cite{Koepf,PWZ}. In fact, the $q$-Zeilberger's algorithm will generate a recursion of
order $3$ for the right-hand side of \eqref{eq:watson8phi7}.
\end{rmk}

We can now obtain a transformation of terminating very-well-poised $_8\phi_7$
series \cite[(2.10.3)]{GR}:
\begin{cor}For $n\geq 0$, there holds
\begin{align}
&\hskip -2mm  {}_8\phi_7
\left[\begin{array}{c}
a,\, qa^{\frac12},\, -qa^{\frac12},\, b,\, c,\, d,\, e,\, q^{-n} \\
a^{\frac12},\, -a^{\frac12},\, aq/b,\, aq/c,\, aq/d,\, aq/e,\, aq^{n+1}
\end{array};q,\,\frac{a^2q^{n+2}}{bcde}\right]   \nonumber\\[5pt]
& =\frac{(aq,\lambda q/e;q)_n}{(aq/e,\lambda q;q)_n}
{}_8\phi_7
\left[\begin{array}{c}
\lambda,\, q\lambda^{\frac12},\, -q\lambda^{\frac12},\, \lambda b/a,\, \lambda c/a,
\, \lambda d/a,\, e,\, q^{-n} \\
\lambda^{\frac12},\, -\lambda^{\frac12},\, aq/b,\, aq/c,\, aq/d,\, \lambda q/e,\,
\lambda q^{n+1}\end{array};
q,\,\frac{aq^{n+1}}{e}
\right],   \label{eq:vwp-8phi7}
\end{align}
where $\lambda=a^2q/bcd$.
\end{cor}

\pf This is just the applications of \eqref{eq:watson8phi7} on
both sides of \eqref{eq:vwp-8phi7}.   \qed

%%%%%%%%%%%%%%%%%%%%%%%%%%%%%%%%%%%%%%%%%%%%%%%%%%%%%%%%%%%%%%%%%%%%%%%%%%%%%%%%%%%%%%%
\section{Bailey's $_{10}\phi_9$ transformation formula}\label{sec:Bailey}
In this section, we show that
Bailey's $_{10}\phi_9$ transformation formula (see \cite[Appendix (III.28)]{GR}) can also be proved
by this method.
\begin{thm}[Bailey's classical $_{10}\phi_9$ transformation]
For $n\geq 0$, there holds
\begin{align}
&\hskip -2mm  {}_{10}\phi_9
\left[\begin{array}{c}
a,\, qa^{\frac12},\, -qa^{\frac12},\, b,\, c,\, d,\, e,\, f,\,\lambda aq^{n+1}/ef,\,q^{-n}
 \\
a^{\frac12},\,-a^{\frac12},\,aq/b,\,aq/c,\,aq/d,\,aq/e,\,aq/f,\,efq^{-n}/\lambda,\,aq^{n+1}
\end{array};q,\,q
\right]   \nonumber\\[5pt]
& =\frac{(aq,aq/ef,\lambda q/e,\lambda q/f;q)_n}{(aq/e,aq/f,\lambda q/ef,\lambda q;q)_n}
\nonumber\\[5pt]
&\quad{}\times {}_{10}\phi_9\left[\begin{array}{c}
\lambda,\, q\lambda^{\frac12},\, -q\lambda^{\frac12},\,
\lambda b/a,\, \lambda c/a,\, \lambda d/a,\, e,\, f,\, \lambda aq^{n+1}/ef,\,q^{-n}
 \\
\lambda^{\frac12},\,-\lambda^{\frac12},\,aq/b,\,aq/c,\,aq/d,\,
\lambda q/e,\, \lambda q/f,efq^{-n}/a,\,\lambda q^{n+1}
\end{array};q,\,q
\right],
\label{eq:10phi9}
\end{align}
where $\lambda=a^2q/bcd$.
\end{thm}

\pf Note that both sides of \eqref{eq:10phi9} are very-well poised. Let
$$
F_{n,k}(a,b,c,d,e,f,q)
=F_k(a,qa^{\frac12},-qa^{\frac12},b,c,d,e,f,\lambda aq^{n+1}/ef,q^{-n};q,q).
$$
Applying the contiguous relation \eqref{eq:contiguous1} with
$a_r=\lambda aq^{n}/ef$ and $a_{r+1}=q^{-n}$, we have
\begin{align}
F_{n,k}(a,b,c,d,e,f,q)-F_{n-1,k}(a,b,c,d,e,f,q)=\alpha_n F_{n-1,k-1}(aq^2,bq,cq,dq,eq,fq,q),
\label{eq:10phi9-rec1}
\end{align}
where
\begin{align*}
\alpha_n &={}-\frac{(1-b)(1-c)(1-d)(1-e)(1-f)}
{(1-aq/b)(1-aq/c)(1-aq/d)(1-aq/e)(1-aq/f)}
\\[5pt]
&\quad\times \frac{(1-aq)(1-aq^2)(1-ef/\lambda)(1-\lambda aq^{2n}/ef)}
{(1-aq^{n})(1-aq^{n+1})(1-efq^{1-n}/\lambda)(1-efq^{-n}/\lambda)q^{n-1}}.
\end{align*}

Let
\begin{align*}
G_{n,k}(a,b,c,d,e,f,q)
&=\frac{(aq,aq/ef,\lambda q/e,\lambda q/f;q)_n}
{(aq/e,aq/f,\lambda q/ef,\lambda q;q)_n} \\[5pt]
&\quad{}\times F_k(\lambda,q\lambda^{\frac12},-q\lambda^{\frac12},
\lambda b/a,\lambda c/a,\lambda d/a,e,f,\lambda aq^{n+1}/ef,q^{-n};q,q).
\end{align*}
Then we may verify that
\begin{align*}
&\hskip -3mm
G_{n,k}(a,b,c,d,e,f,q)-G_{n-1,k}(a,b,c,d,e,f,q)-
\alpha_n G_{n-1,k-1}(aq^2,bq,cq,dq,eq,fq,q)\\[5pt]
&=H_{n,k}-H_{n,k-1},
\end{align*}
where
\begin{align*}
H_{n,k}
&=\frac{(1-a\lambda q^{2n}/ef)(aq,\lambda q/e,\lambda q/f;q)_{n-1}(aq/ef;q)_n}
{(aq/e,aq/f,\lambda q/ef,\lambda;q)_n} \\[5pt]
&\quad{}\times \frac{(1-\lambda q^{k}/a)
(\lambda b/a,\lambda c/a,\lambda d/a, \lambda aq^{n+1}/ef,q^{1-n};q)_k
(\lambda, e,f;q)_{k+1} }
{(q,aq/b,aq/c,aq/d, \lambda q/e, \lambda q/f;q)_k (efq^{-n}/a,\lambda q^{n};q)_{k+1}}.
\tag*{\qed}
\end{align*}

%%%%%%%%%%%%%%%%%%%%%%%%%%%%%%%%%%%%%%%%%%%%%%%%%%%%%%%%%%%%%%%%%%%%%%%%%%%%%%%%%%
\section{Singh's quadratic transformation}
The following quadratic transformation (see \cite[Appendix (III.21)]{GR}) was first proved by Singh \cite{Singh}.
For a more recent proof, see Askey and Wilson \cite{AW}. We will show that it can
also be proved by induction in the same vein as in the previous sections.
\begin{thm}[Singh's quadratic transformation]There holds
\begin{align}
{}_4\phi_3
\left[\begin{array}{c}
a^2,\ b^2,\ c,\ d\\
abq^{\frac12},\, -abq^{\frac12},\, -cd \end{array};q,\,q\right]
={}_4\phi_3
\left[\begin{array}{c}
a^2,\ b^2,\ c^2,\ d^2\\
a^2b^2q,\, -cd,\, -cdq \end{array};q^2,\,q^2\right],
\label{eq:abcd-q}
\end{align}
provided that both series terminate.
\end{thm}

\pf Let $d=q^{-n}$. Then \eqref{eq:abcd-q} may be written as
\begin{align*}
{}_4\phi_3
\left[\begin{array}{c}
a^2,\ b^2,\ c,\ q^{-n}\\
abq^{\frac12},\, -abq^{\frac12},\, -cq^{-n} \end{array};q,\,q\right]
={}_4\phi_3
\left[\begin{array}{c}
a^2,\,b^2,\,c^2,\, q^{-2n}\\
a^2b^2q,\, -cq^{-n},\, -cq^{1-n} \end{array};q^2,\,q^2\right].
%\label{eq:abc-qn}
\end{align*}
Let
\begin{align*}
F_{n,k}(a,b,c,q)&=\frac{(a^2,b^2,c,q^{-n};q)_k}
{(q,abq^{\frac12}, -abq^{\frac12}, -cq^{-n};q)_k} q^k.
\end{align*}
Applying the relation:
\begin{align*}
1-\frac{(1-xy)(1+cy)}{(1-y)(1+cxy)}
&=-\frac{(1+c)(1-x)y}{(1-y)(1+cxy)},
\end{align*}
with $x=q^{k}$ and $y=q^{-n}$, one sees that
\begin{align}
&\hskip -3mm
F_{n,k}(a,b,c,q)-F_{n-1,k}(a,b,c,q) \nonumber\\[5pt]
&=-\frac{(1-a^2)(1-b^2)(1-c^2)q^{1-n}}{(1-a^2b^2q)(1+cq^{-n})(1+cq^{1-n})}
F_{n-1,k-1}(aq^{\frac12},bq^{\frac12},cq,q).
\label{eq:fnkabc}
\end{align}
Replacing $n$ by $n-1$ in \eqref{eq:fnkabc}, we get
\begin{align}
&\hskip -3mm
F_{n-1,k}(a,b,c,q)-F_{n-2,k}(a,b,c,q) \nonumber\\[5pt]
&=-\frac{(1-a^2)(1-b^2)(1-c^2)q^{2-n}}{(1-a^2b^2q)(1+cq^{1-n})(1+cq^{2-n})}
F_{n-2,k-1}(aq^{\frac12},bq^{\frac12},cq,q).
\label{eq:fnkabc1}
\end{align}
From \eqref{eq:fnkabc} and \eqref{eq:fnkabc1} it follows that
\begin{align}
&q(1+cq^{-n})[F_{n,k}(a,b,c,q)-F_{n-1,k}(a,b,c,q)] \nonumber\\[5pt]
&{}-(1+cq^{2-n})[F_{n-1,k}(a,b,c,q)-F_{n-2,k}(a,b,c,q)]\nonumber\\[5pt]
&=-\frac{(1-a^2)(1-b^2)(1-c^2)q^{2-n}}{(1-a^2b^2q)(1+cq^{1-n})}
[F_{n-1,k-1}(aq^{\frac12},bq^{\frac12},cq,q)-F_{n-2,k-1}(aq^{\frac12},bq^{\frac12},cq,q)].
\label{eq:fnkabc2}
\end{align}

If we apply \eqref{eq:fnkabc1} with $a,b,c,k$ replaced by $aq^{\frac12},bq^{\frac12},cq,k-1$,
respectively, then \eqref{eq:fnkabc2},
after dividing both sides by $q(1+cq^{-n})$, may be written as
\begin{align*}
F_{n,k}(a,b,c,q)-\alpha_n F_{n-1,k}(a,b,c,q)+\beta_n F_{n-2,k}(a,b,c,q)
=\gamma_n F_{n-2,k-2}(aq,bq,cq^2,q),
\end{align*}
where
\begin{align*}
\alpha_n &=\frac{(1+q)(1+cq^{1-n})}{q(1+cq^{-n})}, \\[5pt]
\beta_n  &=\frac{1+cq^{2-n}}{q(1+cq^{-n})}, \\[5pt]
\gamma_n &=\frac{(1-a^2)(1-b^2)(1-c^2)(1-a^2q)(1-b^2q)(1-c^2q^2)q^{3-2n}}
{(1-a^2b^2q)(1-a^2b^2q^3)(1+cq^{-n})(1+cq^{1-n})(1+cq^{2-n})(1+cq^{3-n})}.
\end{align*}

Let
\begin{align*}
G_{n,k}(a,b,c,q)&=\frac{(a^2,b^2,c^2,q^{-2n};q^2)_k}
{(q^2;a^2b^2q;q^2)_k (-cq^{-n};q)_{2k}} q^{2k}.
\end{align*}
Then it is easy to verify that
\begin{align*}
&\hskip -3mm
G_{n,k}(a,b,c,q)-\alpha_n G_{n-1,k}(a,b,c,q)+\beta_n G_{n-2,k}(a,b,c,q)
-\gamma_n G_{n-2,k-2}(aq,bq,cq^2,q) \\[5pt]
&=H_{n,k}-H_{n,k-1},
\end{align*}
where
\begin{align*}
H_{n,k}
&=-\frac{(1-q^{2k-1})(a^2,b^2;q^2)_k (q^{4-2n};q^2)_{k-1}(c^2;q^2)_{k+1}q^{2-2n}}
{(q^2;q^2)_{k-1}(a^2b^2q;q^2)_k (-cq^{-n};q)_{2k+2}}. \tag*{\qed}
\end{align*}

%%%%%%%%%%%%%%%%%%%%%%%%%%%%%%%%%%%%%%%%%%%%%%%%%%%%%%%%%%%%%%%%%%%%%%%%%%%%%%%%%%
\section{Schlosser's $C_r$ extension of Jackson's $_8\phi_7$ summation formula}
\label{sec:Schlosser}
Based on a determinant formula of Krattenthaler \cite[Lemma 34]{Krattenthaler},
Schlosser \cite{Schlosser00,Schlosser03} established a $C_r$ extension of Jackson's
$_8\phi_7$ summation formula.
%A generzation of Schlosser's summation formula
%to elliptic hypergeometric series was obtained by Warnaar
%\cite[Theorem 5.1]{Warnaar}.
\begin{thm}[Schlosser's $C_r$ Jackson's sum]\label{thm:Schlosser}
Let $x_1,\ldots,x_r$, $a$, $b$, $c$ and $d$
be indeterminates and let $n$ be a nonnegative integer. Suppose that none of
the denominators in \eqref{eq:Schlosser} vanish. Then
\begin{align}
&\sum_{k_1,\ldots,k_r=0}^{n}\prod_{1\leq i<j\leq r}
\frac{(x_iq^{k_i}-x_jq^{k_j})(1-ax_ix_jq^{k_i+k_j})}{(x_i-x_j)(1-ax_ix_j)}
\prod_{i=1}^{r}\frac{1-ax_i^2 q^{2k_i}}{1-ax_i^2} \nonumber\\[5pt]
&\quad{}\times \prod_{i=1}^{r}
\frac{(ax_i^2,bx_i,cx_i,dx_i,a^2x_iq^{n-r+2}/bcd,q^{-n};q)_{k_i}q^{k_i}}
{(q,ax_iq/b,ax_iq/c,ax_iq/d,bcdx_iq^{r-n-1}/a,ax_i^2q^{n+1};q)_{k_i}}
\nonumber\\[5pt]
&\qquad =\prod_{1\leq i<j\leq r}\frac{1-ax_ix_j q^n}{1-ax_ix_j}
\prod_{i=1}^{r}\frac{(ax_i^2q,aq^{2-i}/bc,aq^{2-i}/bd,aq^{2-i}/cd;q)_n}
{(aq^{2-r}/bcdx_i,ax_iq/b,ax_iq/c,ax_iq/d;q)_n}. \label{eq:Schlosser}
\end{align}
\end{thm}

We shall give an inductive proof of Schlosser's $C_r$ extension of
Jackson's $_8\phi_7$ summation formula. We first give a simple
proof of the $n=1$ case of \eqref{eq:Schlosser},
which we state as the following lemma.
\begin{lem}For $r\geq 0$, there holds
\begin{align}
&\sum_{s_1,\ldots,s_r=0}^{1}\prod_{1\leq i<j\leq r}
\frac{(x_iq^{s_i}-x_jq^{s_j})(1-ax_ix_jq^{s_i+s_j})}
{(x_i-x_j)(1-ax_ix_jq)}
\prod_{i=1}^{r}
\frac{(-1)^{s_i}(bx_i,cx_i,dx_i,a^2x_iq^{3-r}/bcd;q)_{s_i} }
{(ax_iq/b,ax_iq/c,ax_iq/d,bcdx_iq^{r-2}/a;q)_{s_i}} \nonumber\\[5pt]
&\qquad =\prod_{i=1}^{r}\frac{(ax_i^2q,aq^{2-i}/bc,aq^{2-i}/bd,aq^{2-i}/cd;q)_1}
{(aq^{2-r}/bcdx_i,ax_iq/b,ax_iq/c,ax_iq/d;q)_1}. \label{eq:Schlossern=1}
\end{align}
\end{lem}

\pf Multiplying both sides by $\prod_{1\leq i<j\leq r}(x_i-x_j)(1-ax_ix_jq)$,
Equation \eqref{eq:Schlossern=1} may be written as
\begin{align}
&\sum_{s_1,\ldots,s_r=0}^{1}\prod_{1\leq i<j\leq r}
(x_iq^{s_i}-x_jq^{s_j})(1-ax_ix_jq^{s_i+s_j}) \nonumber\\[5pt]
&\quad{}\times \prod_{i=1}^{r}
(-1)^{s_i}(bx_i,cx_i,dx_i,a^2x_iq^{3-r}/bcd;q)_{s_i}
(ax_iq/b,ax_iq/c,ax_iq/d,bcdx_iq^{r-2}/a;q)_{1-s_i} \nonumber\\[5pt]
&\qquad =\prod_{1\leq i<j\leq r}(x_i-x_j)(1-ax_ix_jq)
\prod_{i=1}^{r}
\frac{-bcdx_i(ax_i^2q,aq^{2-i}/bc,aq^{2-i}/bd,aq^{2-i}/cd;q)_1}
{aq^{2-r}}. \label{eq:Schlon=1}
\end{align}
Denote the left-hand side of \eqref{eq:Schlon=1} by $L$.
If $x_i=x_j$ or $x_i=0$ for some $1\leq i,j\leq r$ ($i\neq j$),
then it is easily seen that $L$ is equal to $0$.

If $ax_i x_jq=1$ for some $i\neq j$, then for $0\leq s_i,s_j\leq 1$,
we have
\begin{align*}
1-ax_ix_jq^{s_i+s_j}&=1-q^{s_i+s_j-1},\\[5pt]
(bx_i,cx_i,dx_i,a^2x_iq^{3-r}/bcd;q)_{s_i}
&=\frac{(bx_i,cx_i,dx_i,bcdx_ix_j^2q^{r-1};q)_{s_i}}
{(-bcdx_ix_j^2q^{r-1})^{s_i}},
\\[5pt]
(ax_iq/b,ax_iq/c,ax_iq/d,bcdx_iq^{r-2}/a;q)_{1-s_i}
&=\frac{(bx_j,cx_j,dx_j,bcdx_i^2x_jq^{r-1};q)_{1-s_i}}{
(-bcdx_j^3)^{1-s_i}}.
\end{align*}
It follows that each term on the left-hand side of \eqref{eq:Schlon=1}
subject to $s_i+s_j=1$ is equal to $0$.
Besides, any two terms cancel each other if they have
the same $s_k$ except for $s_i=s_j=0$ and $s_i=s_j=1$, respectively.
Therefore, $L$ is equal to $0$.

If $ax_i^2q=1$ for some $1\leq i\leq r$, then we observe that
any two different terms on the left-hand side of \eqref{eq:Schlon=1} with the same
$s_k$ ($k\neq i$) cancel each other, and hence $L$ is equal to $0$.

Summarizing the above cases, we see that $L$
is divisible by
$$
\prod_{1\leq i<j\leq r}(x_i-x_j)(1-ax_ix_jq) \prod_{i=1}^{r} x_i(1-ax_i^2q).
$$
Now we consider the left-hand side of \eqref{eq:Schlon=1} as a polynomial
in $x_r$. It is easy to see that the coefficient of $x_r^{2r+2}$ is given by
\begin{align*}
&\sum_{s_1,\ldots,s_r=0}^{1}
\left(\prod_{1\leq i<j\leq r-1}(x_iq^{s_i}-x_jq^{s_j})
\prod_{i=1}^{r-1}(ax_iq^{s_i})\right)
(-1)^{s_r}q^{(2r-2)s_r} (a^2q^{3-r})^{s_r} (a^2q^{r+1})^{1-s_r}\\[5pt]
&\quad{}\times \prod_{i=1}^{r-1}
(-1)^{s_i}(bx_i,cx_i,dx_i,a^2x_iq^{3-r}/bcd;q)_{s_i}
(ax_iq/b,ax_iq/c,ax_iq/d,bcdx_iq^{r-2}/a;q)_{1-s_i} \\[5pt]
&\quad{}=0.
\end{align*}
Therefore, if we write
$$
L=P\prod_{1\leq i<j\leq r}(x_i-x_j)(1-ax_ix_jq) \prod_{i=1}^{r} x_i(1-ax_i^2q),
$$
then $P$ is independent of $x_r$. By symmetry, one sees that
$P$ is independent of all $x_i$ ($1\leq i\leq r$).
Now, taking $x_i=q^{i-1}$ ($1\leq i\leq r$),
then
$$
\prod_{1\leq i<j\leq r}(x_iq^{s_i}-x_jq^{s_j})
=\prod_{1\leq i<j\leq r}(q^{i+s_i-1}-q^{j+s_j-1}),
$$
which is equal to $0$ unless $s_1\leq s_2\leq\cdots\leq s_r$. It follows that
\begin{align*}
P&=\sum_{0\leq s_1\leq \cdots\leq s_r \leq 1}\prod_{1\leq i<j\leq r}
\frac{(q^{i+s_i-1}-q^{j+s_j-1})(1-aq^{i+j+s_i+s_j-2})}
{(q^{i-1}-q^{j-1})(1-aq^{i+j-1})}
\prod_{i=1}^{r}\frac{(-1)^{s_i}}{q^{i-1}(1-aq^{2i-1})} \nonumber\\[5pt]
&\quad{}\times \prod_{i=1}^{r}
(bq^{i-1},cq^{i-1},dq^{i-1},a^2q^{i+2-r}/bcd;q)_{s_i}
(aq^i/b,aq^i/c,aq^i/d,bcdq^{i+r-3}/a;q)_{1-s_i} \nonumber\\[5pt]
&=\sum_{k=0}^{r}(-1)^{r-k}q^{\frac{k(k+1-2r)}{2}}
{r\brack k}\frac{1-aq^{2k}}{(aq^k;q)_{r+1}}
(bq^{k},cq^{k},dq^{k},a^2q^{3-r+k}/bcd;q)_{r-k} \nonumber\\[5pt]
&\quad{}\times
(aq/b,aq/c,aq/d,bcdq^{r-2}/a;q)_{k},
\end{align*}
where we have assumed $s_{k}=0$ and $s_{k+1}=1$. It remains to show that
$$
P=(-bcdq^{r-2}/a)^r (aq^{2-r}/bc,aq^{2-r}/bd,aq^{2-r}/cd;q)_r,
$$
or
\begin{align}
&\hskip -3mm
\sum_{k=0}^{r}(-1)^{k}q^{\frac{k(k+1-2r)}{2}}
{r\brack k}\frac{1-aq^{2k}}{(aq^k;q)_{r+1}}
\frac{(aq/b,aq/c,aq/d,bcdq^{r-2}/a;q)_{k}}
{(b,c,d,a^2q^{3-r}/bcd;q)_k}  \nonumber\\[5pt]
&=\frac{(bcdq^{r-2}/a)^r (aq^{2-r}/bc,aq^{2-r}/bd,aq^{2-r}/cd;q)_r}
{(b,c,d,a^2q^{3-r}/bcd;q)_r}. \label{eq:Sch8phi7}
\end{align}
But \eqref{eq:Sch8phi7} follows easily from Jackson's $_8\phi_7$ summation formula
\eqref{eq:8phi7} with parameter substitutions $b\mapsto aq/b$, $c\mapsto aq/c$,
$d\mapsto aq/d$, and $n\mapsto r$.
This proves \eqref{eq:Schlossern=1}. \qed

\noindent{\it Proof of Theorem \ref{thm:Schlosser}. }
Suppose \eqref{eq:Schlosser} holds for $n$.
It is easy to verify that
\begin{align*}
&\frac{(1-a^2x_i q^{n+k_i-r+2}/bcd)(1-bcdx_iq^{k_i+r-n-2})}
{(1-q^{-n+k_i-1})(1-ax_i^2q^{n+k_i+1})} \\[5pt]
&=-\frac{q^{n+1}(1-a^2x_iq^{n-r+2}/bcd)(1-bcdx_iq^{r-n-2}/a)}
{(1-q^{n+1})(1-ax_i^2q^{n+1})} \\[5pt]
&\quad{}+\frac{(1-a^2x_iq^{2n-r+3}/bcd)(1-bcdx_iq^{r-1}/a)
(1-ax_i^2 q^{k_i})(1-q^{k_i})}
{(1-q^{n+1})(1-ax_i^2q^{n+1})(1-q^{-n+k_i-1})(1-ax_i^2q^{n+k_i+1})}.
\end{align*}
Hence,
\begin{align*}
\prod_{i=1}^{r}\frac{(1-a^2x_i q^{n+k_i-r+2}/bcd)(1-bcdx_iq^{k_i+r-n-2})}
{(1-q^{-n+k_i-1})(1-ax_i^2q^{n+k_i+1})}
=\sum_{s_1,\ldots,s_r=0}^{1}
\frac{\alpha_{s_1,\ldots,s_r} (1-ax_i^2 q^{k_i})^{s_i}(1-q^{k_i})^{s_i}}
{(1-q^{-n+k_i-1})^{s_i}(1-ax_i^2q^{n+k_i+1})^{s_i}},
\end{align*}
where
\begin{align*}
\alpha_{s_1,\ldots,s_r}
&=\frac{1}{(1-q^{n+1})^r}
\prod_{i=1}^{r}\frac{(-q^{n+1})^{1-s_i}}{1-ax_i^2q^{n+1}}
(1-a^2x_iq^{n-r+2}/bcd)^{1-s_i}\\[5pt]
&\quad{}\times(1-bcdx_iq^{r-n-2}/a)^{1-s_i}
(1-a^2x_iq^{2n-r+3}/bcd)^{s_i}(1-bcdx_iq^{r-1}/a)^{s_i}.
\end{align*}
It follows that
\begin{align}
&\prod_{1\leq i<j\leq r}
\frac{(x_iq^{k_i}-x_jq^{k_j})(1-ax_ix_jq^{k_i+k_j})}{(x_i-x_j)(1-ax_ix_j)}
\prod_{i=1}^{r}\frac{1-ax_i^2 q^{2k_i}}{1-ax_i^2} \nonumber\\[5pt]
&\quad{}\times \prod_{i=1}^{r}
\frac{(ax_i^2,bx_i,cx_i,dx_i,a^2x_iq^{n-r+3}/bcd,q^{-n-1};q)_{k_i}q^{k_i}}
{(q,ax_iq/b,ax_iq/c,ax_iq/d,bcdx_iq^{r-n-2}/a,ax_i^2q^{n+2};q)_{k_i}}
\nonumber\\[5pt]
&=\sum_{s_1,\ldots,s_r=0}^{1}\beta_{s_1,\ldots,s_r}\prod_{1\leq i<j\leq r}
\frac{(x_iq^{k_i}-x_jq^{k_j})(1-ax_ix_jq^{k_i+k_j})}
{(x_iq^{s_i}-x_jq^{s_j})(1-ax_ix_jq^{s_i+s_j})}
\prod_{i=1}^{r}\frac{1-ax_i^2 q^{2k_i}}{1-ax_i^2q^{2s_i}} \nonumber\\[5pt]
&\quad{}\times \prod_{i=1}^{r}
\frac{(ax_i^2q^{2s_i},bx_iq^{s_i},cx_iq^{s_i},dx_iq^{s_i},
a^2x_iq^{n-r+s_i+2}/bcd,q^{-n};q)_{k_i-s_i}q^{k_i-s_i}}
{(q,ax_iq^{s_i+1}/b,ax_iq^{s_i+1}/c,ax_iq^{s_i+1}/d,
bcdx_iq^{s_i+r-n-1}/a,ax_i^2q^{n+2s_i+1};q)_{k_i-s_i}},  \label{eq:Schlo-ind}
\end{align}
where
\begin{align*}
\beta_{s_1,\ldots,s_r}&
=\alpha_{s_1,\ldots,s_r}\prod_{1\leq i<j\leq r}
\frac{(x_iq^{s_i}-x_jq^{s_j})(1-ax_ix_jq^{s_i+s_j})}{(x_i-x_j)(1-ax_ix_j)}
\prod_{i=1}^{r}\frac{1-ax_i^2 q^{2s_i}}{1-ax_i^2}  \\[5pt]
&\quad{}\times
\frac{(ax_i^2;q)_{2s_i}
(bx_i,cx_i,dx_i;q)_{s_i}q^{s_i}(1-ax_i^2q^{n+s_i+1})^{1-2s_i}(1-q^{-n-1})}
{(ax_iq/b,ax_iq/c,ax_iq/d;q)_{s_i}(bcdx_iq^{r-n-2}/a;q)_{s_i+1}
(a^2x_iq^{n-r+2}/bcd;q)_{1-s_i}}\\[5pt]
&=\prod_{1\leq i<j\leq r}
\frac{(x_iq^{s_i}-x_jq^{s_j})(1-ax_ix_jq^{s_i+s_j})}{(x_i-x_j)(1-ax_ix_j)}
\\[5pt]
&\quad{}\times \prod_{i=1}^{r}
\frac{(-1)^{s_i}(ax_i^2q;q)_{2s_i}
(bx_i,cx_i,dx_i,bcdx_iq^{r-1}/a,a^2x_iq^{2n-r+3}/bcd,q)_{s_i}}
{q^{ns_i}(ax_i^2q^{n+1},bcdx_iq^{r-n-2}/a;q)_{2s_i}
(ax_iq/b,ax_iq/c,ax_iq/d;q)_{s_i}}.
\end{align*}
By the induction hypothesis, the right-hand side of \eqref{eq:Schlo-ind}
is equal to
\begin{align}
&\sum_{s_1,\ldots,s_r=0}^{1}\beta_{s_1,\ldots,s_r}
\prod_{1\leq i<j\leq r}\frac{1-ax_ix_j q^{n+s_i+s_j}}{1-ax_ix_jq^{s_i+s_j}}
\nonumber\\[5pt]
&\quad{}\times\prod_{i=1}^{r}
\frac{(ax_i^2q^{2s_i+1},aq^{2-i}/bc,aq^{2-i}/bd,aq^{2-i}/cd;q)_n}
{(aq^{2-r-s_i}/bcdx_i,ax_iq^{s_i+1}/b,ax_iq^{s_i+1}/c,ax_iq^{s_i+1}/d;q)_n}
\nonumber\\[5pt]
&=\sum_{s_1,\ldots,s_r=0}^{1}\prod_{1\leq i<j\leq r}
\frac{(x_iq^{s_i}-x_jq^{s_j})(1-ax_ix_jq^{n+s_i+s_j})}{(x_i-x_j)(1-ax_ix_j)}
\nonumber\\[5pt]
&\quad{}\times\prod_{i=1}^{r}
\frac{(bx_i,cx_i,dx_i,a^2x_iq^{2n-r+3}/bcd;q)_{s_i}
(ax_i^2q,aq^{2-i}/bc,aq^{2-i}/bd,aq^{2-i}/cd;q)_n}
{(bcdx_iq^{r-n-2})^{s_i}(aq^{2-r}/bcdx_i,ax_iq/b,ax_iq/c,ax_iq/d;q)_{n+s_i}}.
\label{eq:Schlo-end}
\end{align}

Replacing $a$ by $aq^{n}$ in \eqref{eq:Schlossern=1}, we have
\begin{align*}
&\sum_{s_1,\ldots,s_r=0}^{1}\prod_{1\leq i<j\leq r}
\frac{(x_iq^{s_i}-x_jq^{s_j})(1-ax_ix_jq^{n+s_i+s_j})}{(x_i-x_j)(1-ax_ix_j)}
\nonumber\\[5pt]
&\quad{}\times \prod_{i=1}^{r}
\frac{(bx_i,cx_i,dx_i,a^2x_iq^{3-r}/bcd;q)_{s_i} (-1)^{s_i}}
{(ax_iq^{n+1}/b,ax_iq^{n+1}/c,ax_iq^{n+1}/d,bcdx_iq^{r-n-2}/a;q)_{s_i}}
\nonumber\\[5pt]
&\qquad =\prod_{1\leq i<j\leq r}\frac{1-ax_ix_j q^{n+1}}{1-ax_ix_j}
\prod_{i=1}^{r}\frac{(ax_i^2q^{n+1},aq^{n+2-i}/bc,aq^{n+2-i}/bd,aq^{n+2-i}/cd;q)_1 }
{(aq^{n+2-r}/bcdx_i,ax_iq^{n+1}/b,ax_iq^{n+1}/c,ax_iq^{n+1}/d;q)_1}.
\end{align*}
Therefore, the right-hand side of \eqref{eq:Schlo-end} is equal to
\begin{align*}
\prod_{1\leq i<j\leq r}
\frac{1-ax_ix_jq^{n+1}}{1-ax_ix_j}
\prod_{i=1}^{r}
\frac{(ax_i^2q,aq^{2-i}/bc,aq^{2-i}/bd,aq^{2-i}/cd;q)_{n+1}}
{(aq^{2-r}/bcdx_i,ax_iq/b,ax_iq/c,ax_iq/d;q)_{n+1}}.
\end{align*}
Namely, formula \eqref{eq:Schlosser} holds for $n+1$.
This completes the inductive step, and we conclude that \eqref{eq:Schlosser}
holds for all integers $n\geq 0$.
\qed

We end this section with two new allied identities.
\begin{prop}For $n\geq 0$, we have
\begin{align}
&\sum_{s_1,\ldots,s_r=0}^{n}
\prod_{1\leq i<j\leq r}
\frac{(x_iq^{s_i}-x_jq^{s_j})(1-ax_ix_jq^{s_i+s_j})}
{(x_i-x_j)(1-ax_ix_jq^n)}
\prod_{i=1}^{r}\frac{1}{q^{(r-1)s_i}}
=\frac{(n+1)(q^{n+1};q^{n+1})_{r-1}}{q^{n{r\choose 2}}(q;q)_{r-1}},
\label{eq:cr001}\\[5pt]
&\sum_{s_1,\ldots,s_r=0}^{n}\prod_{1\leq i<j\leq r}
\frac{(x_iq^{s_i}-x_jq^{s_j})(1-ax_ix_jq^{s_i+s_j})}
{(x_i-x_j)(1-ax_ix_jq^n)}
\prod_{i=1}^{r}\frac{1}{(-q)^{(r-1)s_i}} \nonumber\\[5pt]
&\qquad=\begin{cases}
\displaystyle
\frac{(-q^{n+1};q^{n+1})_{r-1}}{q^{n{r\choose 2}}(-q;q)_{r-1}},
&\text{if $n\equiv 0\pmod 2$,}\\[10pt]
0, &\text{if $n\equiv 1\pmod 2$.}
\label{eq:cr002}
\end{cases}
\end{align}
\end{prop}

\pf Assume $|q|<1$. Similarly to the proof of \eqref{eq:Schlossern=1},
we can show that the left-hand side of \eqref{eq:cr001} is
independent of $x_1,\ldots,x_r$.
In particular, taking $x_i=q^{mi}$ ($1\leq i\leq r$)
and letting $m\to+\infty$, the left-hand side of \eqref{eq:cr001}
then becomes
\begin{align*}
\sum_{s_1,\ldots,s_r=0}^{n}\prod_{i=1}^{r}\frac{1}{q^{(i-1)s_i}}
=\prod_{i=1}^{r}\sum_{s_i=0}^{n}\frac{1}{q^{(i-1)s_i}},
\end{align*}
which is clearly equal to the right-hand side of \eqref{eq:cr001}.
Similarly we can prove \eqref{eq:cr002}.
\qed

%%%%%%%%%%%%%%%%%%%%%%%%%%%%%%%%%%%%%%%%%%%%%%%%%%%%%%%%%%%%%%%%%%%%%%%%%%%%%%%%%%%
\section{A finite form of Lebesgue's identity and Jacobi's triple product identity}
\label{sec:triple}
Using the same vein of \eqref{eq:fnk-ind}, we can also derive some other identities.
Here we will prove a finite form of Lebesgue's identity and Jacobi's triple product identity.
As usual, the {\it $q$-binomial coefficients} are defined by
\begin{align*}
{n\brack k}=\frac{(q;q)_n}{(q;q)_k(q;q)_{n-k}},\quad n,k\in\mathbb{Z}.
\end{align*}
\begin{thm}
[A finite form of Lebesgue's and Jacobi's identities]
\label{thm:lebes-1}
For $n\geq 0$, there holds
\begin{align}
\sum_{k=0}^{n}{n\brack k}\frac{q^{k(k+1)/2}}{(aq^k;q)_{n+1}}
=\frac{(-q;q)_n}{(a;q^2)_{n+1}}. \label{eq:lebes-1}
\end{align}
\end{thm}

\pf Let
\begin{align*}
F_{n,k}(a,q)={n\brack k}\frac{q^{k(k+1)/2}}{(aq^k;q)_{n+1}}.
\end{align*}
Noticing the trivial relation
\begin{align}
{n\brack k}(1-aq^{n})
={n-1\brack k}(1-aq^{n+k})+{n-1\brack k-1}(1-aq^{k})q^{n-k},
\label{eq:bino-rec}
\end{align}
we have
\begin{align*}
F_{n,k}(a,q)=\frac{1}{1-aq^{n}}F_{n-1,k}(a,q)
+\frac{q^n}{1-aq^n}F_{n-1,k-1}(aq^2,q).
\tag*{\qed}
\end{align*}

Letting $n\rightarrow\infty$ in \eqref{eq:lebes-1}, we immediately get
Lebesgue's identity:
$$
\sum_{k=0}^{\infty}\frac{(a;q)_k}{(q;q)_k}q^{k(k+1)/2}
=(aq;q^2)_\infty(-q;q)_\infty.
$$

It is interesting to note that \eqref{eq:lebes-1} is also a finite form of
{\it Jacobi's triple product identity}
\begin{equation}\label{eq:jacobi}
\sum_{k=-\infty}^{\infty}q^{k^2}z^k=(q^2,-q/z,-qz;q^2)_\infty.
\end{equation}
Indeed, performing parameter replacements $n\rightarrow m+n$, $a\rightarrow -zq^{-2m}$
and $k\rightarrow k+m$ in \eqref{eq:lebes-1}, and noticing the following
relation
$$
\frac{(-zq^{-2m};q^2)_{m+n+1}}{(-zq^{k-m};q)_{m+n+1}}q^{m+k+1\choose 2}
=\frac{(-q^2/z;q^2)_m(-z;q^2)_{n+1}}{(-q/z;q)_{m-k}(-z;q)_{n+k+1}}q^{k^2}z^k,
$$
we obtain
\begin{align}
\sum_{k=-m}^{n}{m+n\brack m+k}
\frac{(-q^2/z;q^2)_m(-z;q^2)_{n+1}}{(-q/z;q)_{m-k}(-z;q)_{n+k+1}}q^{k^2}z^k
=(-q;q)_{m+n}. \label{eq:jacobi-1}
\end{align}
Letting $m,n\rightarrow\infty$ in \eqref{eq:jacobi-1} and applying the
relation
$$
\frac{(-q,q,-q/z,-z;q)_\infty }{(-q^2/z,-z;q^2)_\infty}=(q^2,-q/z,-qz;q^2)_\infty,
$$
we are led to \eqref{eq:jacobi}.

\begin{rmk}
Theorem \ref{thm:lebes-1} is in fact the $b=q^{-n}$ case of
the $q$-analogue of Kummer's theorem (see \cite[Appendix (II.9)]{GR}):
\begin{align*}
_2\phi_1
\left[\begin{array}{c}
a,\,b \\
q,\,aq/b,\end{array};q,\,-q/b
\right] =\frac{(aq;q^2)_\infty(aq^2/b^2;q^2)_\infty (-q;q)_\infty}
{(aq/b;q)_\infty (-q/b;q)_\infty}.
\end{align*}
There also exists another finite form of {\it Lebesgue's identity}
as follows:
\begin{equation}\label{eq:lebesgue1}
\sum_{k=0}^{n}\frac{(a;q)_{k}(q^{-2n};q^2)_k }
{(q;q)_k (q^{-2n};q)_{k}}q^{k}
=\frac{(aq;q^2)_n}{(q;q^2)_n},
\end{equation}
which is the special case $b=0$ of the following identity due to
Andrews-Jain (see \cite{Andrews76,Guo,Jain}):
\begin{align*}
\sum_{k=0}^{n}\frac{(a,b;q)_{k}(q^{-2n};q^2)_k q^{k}}
{(q;q)_k (abq;q^2)_{k} (q^{-2n};q)_{k}}
=\frac{(aq,bq;q^2)_n}{(q,abq;q^2)_n}.
\end{align*}
A bijective proof of Lebesgue's identity was obtained by Bessenrodt \cite{Bessenrodt}
using Sylvester's bijection. It would be interesting to find a bijective proof of \eqref{eq:lebes-1}
and \eqref{eq:lebesgue1}.
Alladi and Berkovich \cite{AB} have given different finite analogues of Jacobi's triple product and
Lebesgue's identities.
\end{rmk}

%%%%%%%%%%%%%%%%%%%%%%%%%%%%%%%%%%%%%%%%%%%%%%%%%%%%%%%%%%%%%%%%%%%%%%%%%%%%%%%%%%%%%%%%%%%
\section{A finite form of Watson's quintuple product identity}\label{sec:quintuple}
{\it Watson's quintuple product identity} (see \cite[p.~147]{GR}) states that
\begin{align}
\sum_{k=-\infty}^{\infty}(z^2 q^{2k+1}-1)
z^{3k+1}q^{k(3k+1)/2}=(q,z,q/z;q)_\infty (qz^2,q/z^2;q^2)_\infty. \label{eq:quint}
\end{align}
We will show that \eqref{eq:quint} follows from the following trivial identity.
\begin{thm}[A finite form of Watson's quintuple product identity]
For $n\geq 0$, there holds
\begin{align}
\sum_{k=0}^{n}(1-z^2 q^{2k+1}){n\brack k}
\frac{(zq;q)_n}{(z^2q^{k+1};q)_{n+1}}z^kq^{k^2}=1.
\label{eq:quint-1}
\end{align}
\end{thm}
\medskip
\pf Let
$$
F_{n,k}(z,q)=(1-z^2 q^{2k+1}){n\brack k}
\frac{(zq;q)_n}{(z^2q^{k+1};q)_{n+1}}z^kq^{k^2}.
$$
Replacing $a$ by $z^2q$ in \eqref{eq:bino-rec}, one sees that
\begin{align*}
F_{n,k}(z,q)
=\frac{1-zq^n}{1-z^2q^{n+1}}F_{n-1,k}(z,q)
+\frac{(1-zq)zq^n}{1-z^2q^{n+1}}F_{n-1,k-1}(zq,q).
\tag*{\qed}
\end{align*}

Making the substitutions $n\rightarrow m+n$, $z\rightarrow -zq^{-m}$
and $k\rightarrow m+k$
in \eqref{eq:quint-1}, and noticing the following relation
\begin{align*}
\frac{(-zq^{1-m};q)_{m+n}}{(z^2q^{k-m+1};q)_{m+n+1}}(-zq^{-m})^k q^{(k+m)^2}
=\frac{(-q/z;q)_{m-1}(-z;q)_{n+1}}{(1/z^2;q)_{m-k}(z^2q;q)_{n+k+1}}
z^{3k-1}q^{k(3k+1)/2},
\end{align*}
we obtain
\begin{align}
\sum_{k=-m}^{n}(1-z^2 q^{2k+1}){m+n\brack m+k}
\frac{(-q/z;q)_{m-1}(-z;q)_{n+1}}{(1/z^2;q)_{m-k}(z^2q;q)_{n+k+1}}
z^{3k-1}q^{k(3k+1)/2}=1.  \label{eq:quintmn}
\end{align}
Letting $m,n\rightarrow \infty$ and applying the relation
$$
\frac{(-z;q)_\infty (-q/z;q)_\infty}{(z^2q;q)_\infty (1/z^2;q)_\infty}
=\frac{-z^2}{(qz^2;q^2)_\infty (z;q)_\infty (q/z^2;q)_\infty (q/z;q)_\infty},
$$
we immediately obtain Watson's quintuple product identity \eqref{eq:quint}.

\begin{rmk}
Note that \eqref{eq:quint-1} is the limiting case
$c\to \infty$ of \eqref{eq:6phi5} with $a=z^2q$
and $b=zq$.
Paule \cite{Paule} has proved the $m=n$ case of \eqref{eq:quintmn}.
On the other hand,
Chen et al.~\cite{CCG} obtained the following
 finite form of Watson's quintuple product
identity:
\begin{align}
\sum_{k=0}^{n}(1+z q^k){n\brack k}
\frac{(z;q)_{n+1}}{(z^2q^k;q)_{n+1}}z^kq^{k^2}=1.
\label{eq:quint-ccg}
\end{align}
Finally, by letting $n\to\infty$ and substituting
$z\to z/q$ and $q\to q^2$ in
\eqref{eq:quint-1}, and $z\to zq$ and $q\to q^2$ in \eqref{eq:quint-ccg}, respectively,
we get,
\begin{align}
1+\sum_{k=1}^{\infty}\frac{z^k q^{2k^2-k}(z^2q^2;q^2)_{k-1}(1-z^2q^{4k})}
{(q^2;q^2)_k}=(-zq;q^2)_\infty (z^2q^4;q^4)_\infty,
\label{eq:ab11}
\end{align}
\begin{align}
\sum_{k=0}^{\infty}\frac{z^k q^{2k^2+k}(z^2q^2;q^2)_k(1+zq^{2k+1})}
{(q^2;q^2)_k}=(-zq;q^2)_\infty (z^2q^4;q^4)_\infty.
\label{eq:ab00}
\end{align}
 Alladi and Berkovich~\cite{AB02} have given the partition interpretations of
 \eqref{eq:ab11} and \eqref{eq:ab00}. It would be interesting to find the corresponding
 combinatorial interpretations of \eqref{eq:quint-1} and \eqref{eq:quint-ccg}.
\end{rmk}

\renewcommand{\baselinestretch}{1}

\end{document}